# Empirical likelihood based testing for regression

**Ingrid Van Keilegom**[*]

*Université catholique de Louvain*
*Institute of Statistics*
*Voie du Roman Pays 20*
*1348 Louvain-la-Neuve*
*Belgium*
*e-mail:* ingrid.vankeilegom@uclouvain.be

**César Sánchez Sellero**[†] **and Wenceslao González Manteiga**[†]

*Universidad de Santiago de Compostela*
*Departamento de Estatística e I.O.*
*Facultad de Matemáticas*
*Campus Sur. 15706 Santiago de Compostela*
*Spain*
*e-mail:* csellero@usc.es; wenceslao@usc.es

**Abstract:** Consider a random vector $(X, Y)$ and let $m(x) = E(Y|X = x)$. We are interested in testing $H_0 : m \in \mathcal{M}_{\Theta,\mathcal{G}} = \{\gamma(\cdot, \theta, g) : \theta \in \Theta, g \in \mathcal{G}\}$ for some known function $\gamma$, some compact set $\Theta \subset \mathbb{R}^p$ and some function set $\mathcal{G}$ of real valued functions. Specific examples of this general hypothesis include testing for a parametric regression model, a generalized linear model, a partial linear model, a single index model, but also the selection of explanatory variables can be considered as a special case of this hypothesis.

To test this null hypothesis, we make use of the so-called marked empirical process introduced by [4] and studied by [16] for the particular case of parametric regression, in combination with the modern technique of empirical likelihood theory in order to obtain a powerful testing procedure. The asymptotic validity of the proposed test is established, and its finite sample performance is compared with other existing tests by means of a simulation study.

[*]Financial support from IAP research networks nr. P5/24 and P6/03 of the Belgian government (Belgian Science Policy) is gratefully acknowledged.

[†]Financial support from the Spanish Ministry of Science and Technology (with additional European FEDER support) through project MTM2005-00820.

**Contents**

## 1. Introduction

Assume that the data $(X_i, Y_i)$ $(i = 1, \ldots, n)$ are independent replications of a random vector $(X, Y)$, where $X$ is a $d$-variate vector and $Y$ is one-dimensional. Let $m(x) = E(Y|X = x)$ for $x \in \mathbb{R}^d$. We are interested in testing

$$H_0 : m \in \mathcal{M}_{\Theta,\mathcal{G}} = \{\gamma(\cdot, \theta, g) : \theta \in \Theta, g \in \mathcal{G}\} \qquad (1.1)$$

for some known function $\gamma$, some compact set $\Theta \subset \mathbb{R}^p$ and some function set $\mathcal{G}$ of real valued functions. Special cases of this general null hypothesis include testing for a parametric model (in which case $m(\cdot) \equiv \gamma(\cdot, \theta)$), a generalized linear model ($m(x) = \gamma(\beta^t x, \alpha)$ with $\theta = (\alpha, \beta)$), a partial linear model ($m(x) = z^t\theta + g(w)$ with $x = (w, z)$), but the test procedure can also be used for e.g. the selection of explanatory variables. Other possibilities not included in the paper are testing for the parametric form of the variance function, the comparison of regression curves, etc.

The idea of the test procedure we develop in this paper is to make use of the so-called marked empirical process introduced by [4] and studied by [16], combined with the modern technique of empirical likelihood theory in order to obtain a powerful testing procedure.

A lot of research on goodness-of-fit tests based on empirical process ideas has been carried out in the last ten years. Starting with the already-mentioned paper by [4] and [16], devoted to testing parametric regression models, the theory was continued with the problem of checking generalized linear models (see [20]), and the selection of variables (see [3]), etc. See also [25] for testing, using the empirical process ideas, other models of interest: partial linear models, reduction of the dimension, models with multidimensional response and heteroscedasticity tests. For the calibration of the critical points associated to the test statistics, two approximations were mainly used: one is based on the bootstrap (see [17]) and the other is based on martingale transformations (see [19]). For alternative testing procedures based on smoothing techniques, see [23], and for recent generalizations of the empirical regression process ideas to dependent data, see the papers by [18], [6], [7] and [8].

The empirical likelihood (see for example the book by [14]) is a technique designed to construct a nonparametric likelihood for parameters of interest in a nonparametric or a semiparametric setting with nice properties typical for the parametric likelihood, as for example the Wilks' theorem and the Bartlett correction recently proved by [1].

There has also been some recent interest in goodness-of-fit tests based on the empirical likelihood in the regression context. [9] propose a sieve empirical likelihood test for testing general varying-coefficient regression models. [11] study the properties of the empirical likelihood in the presence of both finite and infinite dimensional nuisance parameters as well as when the data dimension is high. [2], [21], [12] and [5] propose different tests based on empirical likelihood with conditional moment restrictions, including situations with dependent data.

In this paper we present a new procedure for testing different regression models in a unified way, combining ideas of empirical regression processes and empirical likelihood. See the paper by [10] for a different unified vision of the theory for testing regression models.

The paper is organized as follows. In the next section, the unified testing procedure is explained, and the asymptotic limit of the test is given, under a set of primitive conditions. In Section 3 these primitive conditions are verified for a number of particular regression models. Section 4 contains the results of a small simulation study. The proofs are deferred to the Appendix.

## 2. General test procedure

To explain the idea of the proposed testing procedure, we first focus on the simple case where we test for the parametric form of the regression function, i.e. $m(\cdot) = \gamma(\cdot, \theta)$, and where $X$ is assumed to be one-dimensional. Note that the true value $\theta_0$ of $\theta$ satisfies in that case

$$E[I(X \leq x)(Y - \gamma(X, \theta_0))] = 0 \qquad (2.1)$$

for all $x \in \mathbb{R}$. Let $\hat{\theta}$ be an estimator of $\theta_0$ under $H_0$, obtained from e.g. least squares minimization. The idea of [4] and [16] is to define the marked empirical process based on residuals

$$n^{-1/2} \sum_{i=1}^{n} I(X_i \leq x)[Y_i - \gamma(X_i, \hat{\theta})] \qquad (2.2)$$

and to test $H_0$ by constructing Kolmogorov-Smirnov, Cramér-von Mises, smooth and directional tests based on this process. [16] also obtained the limiting distribution of the proposed tests. In this paper we will extend the above process to a general framework including a large number of common regression models, and we will make use of empirical likelihood theory in the construction of the test statistic.

Let's consider this now in more detail. For the general null hypothesis $H_0$ defined in (1.1), an appropriate extension of the above marked empirical process is given by

$$R_n(u) = n^{-1/2} \sum_{i=1}^{n} L(X_i, u, \hat{\theta}, \hat{g})[Y_i - \gamma(X_i, \hat{\theta}, \hat{g})], \tag{2.3}$$

for $u$ belonging to a compact set $U$, for certain estimators $\hat{\theta}$ and $\hat{g}$ depending on the model under $H_0$, and for a certain appropriate weight function $L$ also depending on the particular form of $H_0$. In Section 3 we will consider the process $R_n(u)$ in detail for a number of particular models (parametric, generalized linear, partial linear, etc.), which will give rise to particular forms for the function $L$ and the set $U$. In many of our examples, $U$ will be equal to $R_X$, the support of $X$, which we also assume to be compact.

We will now introduce the empirical likelihood test based on the marked empirical process $R_n(u)$. For any bivariate distribution $F$ defined on the support of $(X, Y)$, define the likelihood $L(F)$ by (we restrict for simplicity to the definition for $d = 1$, the extension to general $d$ being straightfoward, but notationally more heavy)

$$L(F) = \prod_{i=1}^{n}[F(X_i, Y_i) - F(X_i-, Y_i) - F(X_i, Y_i-) + F(X_i-, Y_i-)].$$

Next, define the empirical likelihood ratio by (for any $d \geq 1$)

$$\begin{aligned}
&EL(u, \hat{\theta}, \hat{g}) \\
&= \frac{\sup\{L(F) : E_F[L(X, u, \hat{\theta}, \hat{g})(Y - \gamma(X, \hat{\theta}, \hat{g}))] = 0\}}{\sup\{L(F)\}} \\
&= \frac{\sup\left\{\prod_{i=1}^{n} w_i : w_i \geq 0, \sum_{i=1}^{n} w_i = 1, \sum_{i=1}^{n} w_i L(X_i, u, \hat{\theta}, \hat{g})[Y_i - \gamma(X_i, \hat{\theta}, \hat{g})] = 0\right\}}{\prod_{i=1}^{n} n^{-1}} \\
&= \sup\left\{\prod_{i=1}^{n}(nw_i) : w_i \geq 0, \sum_{i=1}^{n} w_i = 1, \sum_{i=1}^{n} w_i L(X_i, u, \hat{\theta}, \hat{g})[Y_i - \gamma(X_i, \hat{\theta}, \hat{g})] = 0\right\},
\end{aligned} \tag{2.4}$$

where the maximum over the empty set is defined to be zero, and let $\ell(u, \hat{\theta}, \hat{g}) = -2 \log EL(u, \hat{\theta}, \hat{g})$. The proposed test statistics are now given by

$$S_n = \sup_{u \in U} \ell(u, \hat{\theta}, \hat{g}), \qquad T_n = \int_U \ell(u, \hat{\theta}, \hat{g}) \, d\mu_n(u), \tag{2.5}$$

for some measure $\mu_n$ such that

$$\sup_{u \in U} |\mu_n(u) - \mu(u)| = o_P(1),$$

where the measure $\mu$ does not depend on $n$. When $U = R_X$, a natural choice for $\mu_n$ is
$$d\mu_n(\cdot) = w(\cdot)d\hat{F}_X(\cdot),$$
where $w(\cdot)$ is an appropriate weight function and $\hat{F}_X(\cdot)$ is the empirical distribution function of $X$. The measure $\mu(\cdot)$ equals $d\mu(\cdot) = w(\cdot)dF_X(\cdot)$ in that case, with $F_X(\cdot) = P(X \le \cdot)$.

In order to obtain the limiting behavior of the process $\ell(u, \hat{\theta}, \hat{g})$ ($u \in U$), use will be made of the results in [11]. In that paper the authors give primitive conditions under which the empirical likelihood statistic converges (for a fixed $u$). Their result can be easily generalized to the current situation of a process in $u$. We will show this below. For this purpose, let $\ell^\infty(U)$ be the space of all bounded functions from $U$ to $\mathbb{R}$ equiped with the supremum metric.

Let $\theta_0$ be the true value of $\theta \in \Theta$, and $g_0$ be the true function in $\mathcal{G}$. The result gives the asymptotic behavior of the test statistics $S_n$ and $T_n$ both under $H_0$ and under the fixed alternative hypothesis $H_A : m \notin \mathcal{M}_{\Theta, \mathcal{G}}$.

**Theorem 2.1.** *(a) Assume that under $H_0$,*

- *(C0) $P(EL(u, \hat{\theta}, \hat{g}) = 0$ for some $u \in U) \to 0$ as $n \to \infty$.*
- *(C1) The process $R_n(u)$ can be written as*
  $R_n(u) = n^{-1/2} \sum_{i=1}^n Q(X_i, Y_i, u, \theta_0, g_0) + o_P(1)$, *uniformly in $u \in U$, for some function $Q$ defined on $R_X \times \mathbb{R} \times U \times \Theta \times \mathcal{G}$ that satisfies $E[Q(X, Y, u, \theta_0, g_0)] = 0$ for all $u$, and $R_n(u)$ converges weakly to a zero-mean Gaussian process $V(u)$ ($u \in U$).*
- *(C2) $\sup_{u \in U} |n^{-1} \sum_{i=1}^n L^2(X_i, u, \hat{\theta}, \hat{g})(Y_i - \gamma(X_i, \hat{\theta}, \hat{g}))^2 - T(u)| \to_P 0$ for some function $T(u)$ ($u \in U$), such that $\inf_{u \in U} T(u) > 0$ and $\sup_{u \in U} T(u) < \infty$.*
- *(C3) $\sup_{u \in U} \max_{1 \le i \le n} |L(X_i, u, \hat{\theta}, \hat{g})(Y_i - \gamma(X_i, \hat{\theta}, \hat{g}))| = o_P(n^{1/2})$.*

*Then, under $H_0$, the process $\ell(u, \hat{\theta}, \hat{g})$ ($u \in U$) converges weakly in $\ell^\infty(U)$ to $W^2(u)$, where $W(u) = V(u)/\sqrt{T(u)}$, and hence*
$$S_n \to_d \sup_{u \in U} W^2(u), \quad T_n \to_d \int_U W^2(u)\, d\mu(u).$$

*(b) Assume that under $H_A$, condition (C2) holds, and in addition*

- *(C1') $n^{1/2} \left( \inf_{u \in U_0} |R_n(u)| \right)^{-1} = O_P(1)$, for a subset $U_0$ of $U$ of positive $\mu$-measure.*
- *(C3') $\sup_{u \in U} \max_{1 \le i \le n} |L(X_i, u, \hat{\theta}, \hat{g})^2 (Y_i - \gamma(X_i, \hat{\theta}, \hat{g}))^2| = o_P(n^{1/2})$.*

*Then, under $H_A$, for any $c > 0$,*
$$P(S_n > c) \to 1, \quad P(T_n > c) \to 1,$$
*as $n$ tends to infinity.*

Note that condition (C1') is basically saying that for some values of $u$, $R_n(u)$ is biased. Indeed, if $R_n(u)$ is unbiased for all $u$, then $|R_n(u)| = O_P(1)$ for all $u \in U$ and condition (C1') does not hold in that case. Also note that condition (C3') is stronger than (C3). Consider the common case where $L(X, u, \theta, g)$ is a bounded function, and write

$$Y - \gamma(X, \hat{\theta}, \hat{g}) = [Y - \gamma(X, \theta, g)] + [\gamma(X, \theta, g) - \gamma(X, \hat{\theta}, \hat{g})].$$

Then, the second term can be taken care of using standard arguments, and hence it follows from Lemma 11.2 in [14] that (C3') holds true if the fourth moment of $Y - \gamma(X, \theta, g)$ is finite, whereas (C3) only requires the finiteness of the second moment.

Since the limiting distribution of $S_n$ and $T_n$ are rather complicated, we propose to work with the following bootstrap approximation. From Theorem 2.1(a) it follows that the log-likelihood function $\ell(u, \hat{\theta}, \hat{g})$ can be approximated by $R_n(u)^2/T(u)$ and that $R_n(u)$ can be written as

$$R_n(u) = n^{-1/2} \sum_{i=1}^{n} Q(X_i, Y_i, u, \theta_0, g_0) + o_P(1),$$

uniformly in $u \in U$. Define random variables $V_1, \ldots, V_n$ that are independent of the data $(X_i, Y_i)$ ($i = 1, \ldots, n$), and that are independent and identically distributed such that $E(V_i) = 0$, $\text{Var}(V_i) = 1$ and $|V_i| \leq c < \infty$ for some finite $c$. Then, let

$$\tilde{R}_n^*(u) = n^{-1/2} \sum_{i=1}^{n} Q(X_i, Y_i, u, \theta_0, g_0) V_i.$$

From the multiplier central limit theorem (see e.g. [22], p. 179), it follows that $\tilde{R}_n^*(u)$ converges to the same process as the original process $R_n(u)$. Since the process $\tilde{R}_n^*(u)$ contains unknown quantities ($\theta_0$, $g_0$ and possibly also the function $Q$ itself), we cannot use it in our bootstrap procedure, and we therefore define

$$R_n^*(u) = n^{-1/2} \sum_{i=1}^{n} \hat{Q}(X_i, Y_i, u, \hat{\theta}, \hat{g}) V_i,$$

where $\hat{Q}$ is a suitable estimator of the function $Q$, whose definition depends on the particular model at hand. We will give the precise definition of this estimator for each of the models we will consider in Section 3. Now, suppose the following assumption holds true :

(C4) $\sup_{u \in U} |R_n^*(u) - \tilde{R}_n^*(u)| = o_P(1)$.

Then, $R_n^*(u)$ also converges to the same process as the original process $R_n(u)$. This way of approximating the process $R_n(u)$ by means of the bootstrap has been used in a number of other papers, see e.g. [3], [25], among others.

It remains now to define an estimator of the variance $T(u)$. Let

$$\hat{T}(u) = n^{-1} \sum_{i=1}^{n} L(X_i, u, \hat{\theta}, \hat{g})(Y_i - \gamma(X_i, \hat{\theta}, \hat{g}))^2.$$

Since condition (C2) assures that $\sup_u |\hat{T}(u) - T(u)| = o_P(1)$, we have the following result concerning the bootstrap approximation of the process $W^2(u)$, $u \in U$:

**Theorem 2.2.** *Assume (C0)–(C4). Then, under $H_0$ and conditionally on the data $(X_i, Y_i)$ ($1 \leq i \leq n$), the process $R_n^*(u)^2/\hat{T}(u)$ converges weakly in $\ell^\infty(U)$ to $W^2(u)$.*

As a consequence, the critical values of the test statistics $S_n$ and $T_n$ can be approximated by the $(1-\alpha)$-th quantiles of the distributions of

$$S_n^* = \sup_{u \in U} \frac{R_n^*(u)^2}{\hat{T}(u)}$$

and

$$T_n^* = \int_U \frac{R_n^*(u)^2}{\hat{T}(u)} \, d\mu_n(u),$$

respectively, given the data $(X_i, Y_i)$ ($1 \leq i \leq n$).

## 3. Application of general test to specific models

In this section we consider a number of particular models and apply the test procedure developed in Section 2 to each of these models.

### 3.1. Parametric models

Consider the null hypothesis $H_0^P : m(\cdot) \equiv \gamma(\cdot, \theta)$ for some $\theta \in \Theta$, where $X$ is supposed to be one-dimensional ($d = 1$). As we have seen in Section 2, [4] and [16] introduced the marked empirical process defined in (2.2). We focus here on the following slight variation of (2.2):

$$R_n^P(x) = n^{-1/2} \sum_{i=1}^n I(X_i \in J_x)[Y_i - \gamma(X_i, \hat{\theta})],$$

where $x \in R_X$, and

$$J_x = \begin{cases} \{t : t \geq x\} & x \leq a \\ \{t : t \leq x\} & x > a, \end{cases} \tag{3.1}$$

for some $a$ in the interior of $R_X$ (e.g. $a$ is the median of $X$). Note that $R_n^P(x)$ equals $R_n(x)$ with $L(X, x, \theta, g) = I(X \in J_x)$ and with $\gamma(x, \theta, g)$ being a function of $\theta$ only.

The marked empirical process proposed by [16] accumulates data from left to right (i.e. it is based on the indicator $I(X \leq x)$), whereas we replaced $I(X \leq x)$ by $I(X \in J_x)$, in order to obtain more stable results for small values of $x$. Note that when we take $a = \inf\{x : x \in R_X\}$, the left endpoint of the support of $X$, we obtain the process studied by [16].

To show the validity of conditions (C0)–(C3) we will use the results in [16].

Note that [2] also use an empirical likelihood approach to test the appropriateness of a parametric regression model. In contrast to our method, they use a local approach, based on comparing a kernel estimator of the regression function $m(x)$ with an estimator of $m(x)$ under the null model. As a consequence their results depend on the choice of a smoothing parameter.

Let $T(x) = E[\text{Var}(Y|X)I(X \in J_x)]$. The following conditions are needed for the main result. Note that the choice $a = \inf\{x : x \in R_X\}$ (i.e. the process of [16]) is excluded by assumption $(P1)$.

**(P1)** $E(Y^2) < \infty$, $\inf_{x \in R_X} T(x) > 0$ and $\sup_{x \in R_X} T(x) < \infty$.

**(P2)** Under $H_0$, $\hat{\theta}$ admits an expansion

$$\hat{\theta} - \theta_0 = n^{-1} \sum_{i=1}^n h(X_i, Y_i, \theta_0) + o_P(n^{-1/2}),$$

where $E[h(X, Y, \theta_0)] = 0$ and $H(\theta_0) := E[h(X, Y, \theta_0)h^t(X, Y, \theta_0)]$ exists and is positive definite.

**(P3)(i)** For each $x \in R_X$, $\gamma(x, \theta)$ is continuously differentiable at each $\theta$ in int$(\Theta)$, and $\sup_{x \in R_X} |\gamma(x, \theta_0)| < \infty$.
**(ii)** The components of the vector $\frac{\partial}{\partial \theta}\gamma(x, \theta)$ are bounded, uniformly in $x \in R_X$ and $\theta \in \Theta$.

Condition $(P2)$ is satisfied for the least squares estimator and any of its common robust modifications.

Also, define

$$G(x, \theta) = \int_{J_x} \frac{\partial}{\partial \theta}\gamma(u, \theta) \, dF_X(u).$$

**Theorem 3.1.** *Assume $(P1)$–$(P3)$. Then, under $H_0^P$, conditions (C0)–(C3) hold true for*

$$Q(\bar{x}, \bar{y}, x, \theta) = I(\bar{x} \in J_x)[\bar{y} - \gamma(\bar{x}, \theta)] - G^t(x, \theta)h(\bar{x}, \bar{y}, \theta)$$

*and $V(x)$ a zero-mean Gaussian process with covariance function*

$$\begin{aligned}\text{Cov}(V(x_1), V(x_2)) = &\, E[\text{Var}(Y|X)I(X \in J_{x_1} \cap J_{x_2})] + G^t(x_1, \theta_0)H(\theta_0)G(x_2, \theta_0) \\ &- G^t(x_1, \theta_0)E[I(X \in J_{x_2})(Y - \gamma(X, \theta_0))h(X, Y, \theta_0)] \\ &- G^t(x_2, \theta_0)E[I(X \in J_{x_1})(Y - \gamma(X, \theta_0))h(X, Y, \theta_0)].\end{aligned}$$

In order to apply the bootstrap procedure described in Section 2, let $\tilde{R}_n^{P*}(x) = n^{-1/2} \sum_{i=1}^n Q(X_i, Y_i, x, \theta_0) V_i$ and $R_n^{P*}(x) = n^{-1/2} \sum_{i=1}^n \hat{Q}(X_i, Y_i, x, \hat{\theta}) V_i$, where the random variables $V_i$ are as in Section 2,

$$\hat{Q}(\bar{x}, \bar{y}, x, \theta) = I(\bar{x} \in J_x)[\bar{y} - \gamma(\bar{x}, \theta)] - \hat{G}^t(x, \theta)\hat{h}(\bar{x}, \bar{y}, \theta),$$

I. Van Keilegom et al./Empirical likelihood based testing for regression        589$\hat{G}(x,\theta) = n^{-1}\sum_{i=1}^{n} I(X_i \in J_x)\frac{\partial}{\partial\theta}\gamma(X_i,\theta)$, and $\hat{h}$ is an appropriate estimator of the function $h$ satisfying

**(P4)** $n^{-1/2}\sum_{i=1}^{n}[\hat{h}(X_i,Y_i,\hat{\theta}) - h(X_i,Y_i,\theta_0)]V_i = o_P(1)$.

We also need to assume that

**(P5)** For each $x \in R_X$, $\gamma(x,\theta)$ is twice continuously differentiable with respect to the components of $\theta$ at each $\theta$ in int($\Theta$), and all the partial derivatives of order two are bounded, uniformly in $x \in R_X$ and $\theta \in \Theta$.

**Theorem 3.2.** *Assume (P1)–(P5). Then, under $H_0^P$, condition (C4) holds true, and hence the process $R_n^{P*}(x)$, $x \in R_X$ has the same limiting distribution as the process $R_n^P(x)$.*

Note that condition (P4) is satisfied when e.g. $\gamma(x,\theta) = \theta_0 + \theta_1 x$, and $\hat{\theta}$ is the least squares estimator of $\theta$, since in that case (we consider only $\theta_1$ and write $h(x,y,\theta) = (h_0(x,y,\theta), h_1(x,y,\theta))^t$)

$$h_1(x,y,\theta) = \frac{1}{\sigma_X^2}(x-\mu_X)(y-\mu_Y) - \frac{\sigma_{XY}}{\sigma_X^4}(x-\mu_X)^2,$$

where $\mu_X = E(X)$, $\mu_Y = E(Y)$, $\sigma_X^2 = \text{Var}(X)$ and $\sigma_{XY} = \text{Cov}(X,Y)$. Define $\hat{h}(x,y,\theta)$ by replacing $\mu_X$, $\mu_Y$, $\sigma_X^2$ and $\sigma_{XY}$ by their sample versions. Then, the left hand side of condition (P4) can be decomposed in a sum of i.i.d. terms, plus a remainder term. Consider e.g. the following term (the other terms are similar and lead to the same asymptotic order) :

$$\left\{\frac{1}{\hat{\sigma}_X^2} - \frac{1}{\sigma_X^2}\right\}n^{-1/2}\sum_{i=1}^{n}(X_i-\mu_X)(Y_i-\mu_Y)V_i = O_P(n^{-1/2}),$$

since $E[(X_i-\mu_X)(Y_i-\mu_Y)V_i] = E[(X_i-\mu_X)(Y_i-\mu_Y)]E(V_i) = 0$, and since $\hat{\sigma}_X^2 - \sigma_X^2 = O_P(n^{-1/2})$.

Other bootstrap procedures can be used as well (see e.g. [17] for a wild bootstrap procedure).

### 3.2. Generalized linear models

Under this model it is assumed that $H_0^{GLM} : m(\cdot) \equiv \gamma(\beta^t \cdot, \alpha)$ for some $\theta = (\alpha,\beta) \in \Theta \subset \mathbb{R}^{a+d}$, $d$ being the dimension of $X$. [20] proposed the following marked empirical process (except for the indicator $I(\hat{\beta}^t X_i \in J_u)$, which is replaced by $I(\hat{\beta}^t X_i \leq u)$ in their paper)

$$R_n^{GLM}(u) = n^{-1/2}\sum_{i=1}^{n} I(\hat{\beta}^t X_i \in J_u)[Y_i - \gamma(\hat{\beta}^t X_i, \hat{\alpha})],$$

for some estimators $\hat{\alpha}$ and $\hat{\beta}$, which is of the general form given by $R_n(u)$ with $L(X,u,\theta,g)$ corresponding to $I(\beta^t X \in J_u)$ and $U = \{\beta_0^t x : x \in R_X\}$. The set $J_u$ is as in (3.1), with $a$ belonging now to the interior of the set $U$. The motivation for this process comes from the fact that under $H_0^{GLM}$, $E\{I(\beta_0^t X \in J_u)[Y - \gamma(\beta_0^t X, \alpha_0)]\} = 0$ for all $u$.

Let $T(u) = E[\text{Var}(Y|\beta_0^t X)I(\beta_0^t X \in J_u)]$. The following assumptions are needed for showing the validity of (C0)–(C3).

**(GLM1)** $E(Y^2) < \infty$, $\inf_{u \in U} T(u) > 0$ and $\sup_{u \in U} T(u) < \infty$.

**(GLM2)** Under $H_0$, $(\hat{\alpha}, \hat{\beta})$ admits an expansion

$$(\hat{\alpha}^t, \hat{\beta}^t)^t - (\alpha_0^t, \beta_0^t)^t = n^{-1} \sum_{i=1}^n k(X_i, Y_i, \alpha_0, \beta_0) + o_P(n^{-1/2}),$$

where $E[k(X, Y, \alpha_0, \beta_0)] = 0$ and $K(\alpha_0, \beta_0) := E[k(X, Y, \alpha_0, \beta_0)k^t(X, Y, \alpha_0, \beta_0)]$ exists and is positive definite.

**(GLM3)(i)** For each $x \in R_X$, $\gamma(\beta^t x, \alpha)$ is continuously differentiable at each $(\alpha, \beta)$ in $\text{int}(\Theta)$, and $\sup_{x \in R_X} |\gamma(\beta_0^t x, \alpha_0)| < \infty$.
**(ii)** The components of the vector $\frac{\partial}{\partial(\alpha,\beta)} \gamma(\beta^t x, \alpha)$ are bounded, uniformly in $x \in R_X$ and $(\alpha, \beta) \in \Theta$.

**(GLM4)** The function $E[\text{Var}(Y|X)I(\beta_0^t X \in J_u)]$ is uniformly continuous in $u \in U$.

Also, define

$$G(u, \alpha, \beta) = E\Big[\frac{\partial}{\partial(\alpha, \beta)} \gamma(\beta^t X, \alpha) I(\beta^t X \in J_u)\Big].$$

**Theorem 3.3.** *Assume (GLM1)–(GLM4). Then, under $H_0^{GLM}$, conditions (C0)–(C3) hold true for $V(u)$ a zero-mean Gaussian process with covariance function ($\theta_0 = (\alpha_0, \beta_0)$)*

$$\begin{aligned}\text{Cov}(V(u_1), V(u_2)) = {} & E[\text{Var}(Y|\beta_0^t X) I(\beta_0^t X \in J_{u_1} \cap J_{u_2})] + G^t(u_1, \theta_0) K(\theta_0) G(u_2, \theta_0) \\ & - G^t(u_1, \theta_0) E[I(\beta_0^t X \in J_{u_2})\{Y - \gamma(\beta_0^t X, \alpha_0)\} k(X, Y, \theta_0)] \\ & - G^t(u_2, \theta_0) E[I(\beta_0^t X \in J_{u_1})\{Y - \gamma(\beta_0^t X, \alpha_0)\} k(X, Y, \theta_0)].\end{aligned}$$

For showing the validity of the bootstrap approximation, first note that $R_n^{GLM}(u) = n^{-1/2} \sum_{i=1}^n Q(X_i, Y_i, u, \theta_0) + o_P(1)$, uniformly in $u \in U$, where

$$Q(x, y, u, \theta) = I(\beta^t x \in J_u)[y - \gamma(\beta^t x, \alpha)] - G^t(u, \alpha, \beta) k(x, y, \alpha, \beta).$$

Define

$$\hat{Q}(x, y, u, \theta) = I(\beta^t x \in J_u)[y - \gamma(\beta^t x, \alpha)] - \hat{G}^t(u, \alpha, \beta) \hat{k}(x, y, \alpha, \beta),$$

where $\hat{G}(u, \alpha, \beta)$ is the sample version of $G(u, \alpha, \beta)$ and $\hat{k}(x, y, \alpha, \beta)$ satisfies

**(GLM5)** $n^{-1/2} \sum_{i=1}^n [\hat{k}(X_i, Y_i, \hat{\alpha}, \hat{\beta}) - k(X_i, Y_i, \alpha_0, \beta_0)] V_i = o_P(1)$.

Also, assume that

**(GLM6)** For each $x \in R_X$, $\gamma(\beta^t x, \alpha)$ is twice continuously differentiable with respect to the components of $(\alpha, \beta)$ at each $(\alpha, \beta)$ in int$(\Theta)$, and all the partial derivatives of order two are bounded, uniformly in $x \in R_X$ and $(\alpha, \beta) \in \Theta$.

**Theorem 3.4.** *Assume* (GLM1)–(GLM6). *Then, under* $H_0^{GLM}$, *condition (C4) holds true, and hence the process* $R_n^{GLM*}(u) = n^{-1/2} \sum_{i=1}^n \hat{Q}(X_i, Y_i, u, \hat{\theta}) V_i$, *has the same limiting distribution as the process* $R_n^{GLM}(u)$, $u \in U$.

### 3.3. Selection of variables

For the purpose of this subsection we write $X = (X_1, \ldots, X_d) = (W, Z)$, where $W$, respectively $Z$, is a random vector of dimension $d_w$, respectively $d_z$, and $d = d_w + d_z$. Also, write $R_X = R_W \times R_Z$. Consider the null hypothesis $H_0^{SV}$ : $E(Y|X) = g(W)$, where $g(w) = E(Y|W = w)$, i.e. the response depends only on the covariate vector $W$ and not on $Z$. [3] considered this testing problem and used the fact that under $H_0^{SV}$, $E\{I(X \in J_x) f_W(W)[Y - g(W)]\} = 0$ for all $x \in R_X$, where $f_W$ is the probability density of $W$, to propose the following test statistic :

$$R_n^{SV}(x) = n^{-1/2} \sum_{i=1}^n I(X_i \in J_x) \hat{f}_W(W_i)[Y_i - \hat{g}(W_i)],$$

where $J_x = \prod_{j=1}^d J_{x_j}$, each $J_{x_j}$ is as in (3.1) for some $a_j$ in the interior of the support of $X_j$, $I(X_i \in J_x) = I(X_{i1} \in J_{x_1}, \ldots, X_{id} \in J_{x_d})$,

$$\hat{f}_W(w) = (nh^{d_w})^{-1} \sum_{i=1}^n K\Big(\frac{W_i - w}{h}\Big), \qquad (3.2)$$

$$\hat{g}(w) = \hat{f}_W(w)^{-1}(nh^{d_w})^{-1} \sum_{i=1}^n K\Big(\frac{W_i - w}{h}\Big) Y_i, \qquad (3.3)$$

$K(u) = \prod_{j=1}^{d_w} k(u_j)$, $k$ is a univariate kernel and $h$ a bandwidth. Note that as before, we have replaced the indicator $I(X_i \leq x)$ by $I(X_i \in J_x)$, but the fundamental idea of [3] remains the same.

Define $T(x) = E[\text{Var}(Y|X) f_W^2(W) I(X \in J_x)]$ and $q(z|w) = P(Z \in J_z|W = w)$, and consider the following assumptions :

**(SV1)** $f_W \in \mathcal{F}_\lambda^\infty$, $g \in \mathcal{F}_\tau^2$ and $q(z|\cdot) \in \mathcal{F}_v^\infty$ for all $z \in R_Z$ and for some $\lambda, \tau, v > 0$.

**(SV2)** $k : \mathbb{R} \to \mathbb{R}$ is an even function of uniformly bounded variation, satisfying $k(u) = O((1 + |u|^{\alpha+1+\varepsilon})^{-1})$ for some $\varepsilon > 0$, $\int u^i k(u)\, du = \delta_{i0}$, for $i = 0, \ldots, \alpha-1$, where $\delta_{ij}$ is Kronecker's delta and $\alpha = \ell + t - 1$, with $\ell - 1 < \lambda \leq \ell$ and $t - 1 < \tau \leq t$.

**(SV3)** $(nh^{d_w})^{-1} + nh^{2\min(\tau,\lambda+1)} \to 0$ as $n \to \infty$.

**(SV4)** $E\{|Y - g(W)|^{2+\delta}\} < \infty$ for some $\delta > 0$.

**(SV5)** $\inf_{x \in R_X} T(x) > 0$ and $\sup_{x \in R_X} T(x) < \infty$.

Assumptions $(SV1)$–$(SV4)$ are taken from [3], and ensure the validity of condition (C1). The class $\mathcal{F}_\beta^\alpha$ (for $\alpha, \beta > 0$) defines a class of smooth functions, satisfying certain moment conditions. We refer to [3], p. 1474, for its formal definition.

**Theorem 3.5.** *Assume $(SV1)$–$(SV5)$. Then, under $H_0^{SV}$, conditions (C0)–(C3) hold true for $V(x)$ a zero-mean Gaussian process with covariance function (where $x_j = (w_j, z_j)$, $j = 1, 2$)*

$$Cov(V(x_1), V(x_2)) = E\Big\{(Y - g(W))^2 f_W^2(W) I(W \in J_{w_1} \cap J_{w_2})$$
$$\times [I(Z \in J_{z_1}) - q(z_1|W)][I(Z \in J_{z_2}) - q(z_2|W)]\Big\}.$$

For the bootstrap approximation, [3] showed that $R_n^{SV}(x)$ can be written as

$$R_n^{SV}(x) = n^{-1/2} \sum_{i=1}^n \Big\{ I(X_i \in J_x) f_W(W_i)(Y_i - g(W_i)) -$$
$$- r(x, W_i) f_W(W_i)(Y_i - g(W_i)) \Big\} + o_P(1),$$

where $r(x, W_i) = I(W_i \in J_w) q(z|W_i)$, and they proposed to estimate $r(x, W_i)$ by

$$\hat{r}(x, W_i) = \frac{1}{nh^{d_w} \hat{f}_W(W_i)} \sum_{j=1}^n K\Big(\frac{W_i - W_j}{h}\Big) I(X_j \in J_x).$$

The consistency of the bootstrap approximation obtained by multiplying each term in the above representation by $V_i$, and by replacing the unknown functions in this representation by their corresponding estimators, has been established by the same authors. We refer to their paper for details about the proof and the assumptions under which this approximation is valid.

### 3.4. Partial linear models

We continue to write $X = (X_1, \ldots, X_d) = (W, Z)$ as in the previous subsection, and consider now the following null hypothesis :

$$H_0^{PL} : E(Y|X) = \theta_0^t Z + g(W),$$

i.e. $Y = \theta_0^t Z + g(W) + \varepsilon$, where $E(\varepsilon|W, Z) = 0$, and where the function $g$ is completely unknown. This problem has been studied by [3] and [25]. Note that under $H_0^{PL}$, $E\{Y - m(W) - \theta_0^t[Z - m_Z(W)]|X\} = 0$, where $m(W) = E(Y|W) =$

*I. Van Keilegom et al./Empirical likelihood based testing for regression* 593......

$\theta_0^t E(Z|W) + g(W)$ and $m_Z(W) = E(Z|W)$. This suggests to define the following marked empirical process :

$$R_n^{PL}(x) = n^{-1/2} \sum_{i=1}^n I(X_i \in J_x) w(W_i) \Big\{ Y_i - \hat{m}(W_i) - \hat{\theta}^t [Z_i - \hat{m}_Z(W_i)] \Big\},$$

$x \in R_X$, where $w(\cdot)$ is an appropriate weight function (e.g. $w(W) = f_W(W)$),

$$\hat{m}_Z(w) = \hat{f}_W(w)^{-1}(nh^{d_w})^{-1} \sum_{i=1}^n K\Big(\frac{W_i - w}{h}\Big) Z_i,$$

$K(u) = \prod_{j=1}^{d_w} k(u_j)$, $k$ is a univariate kernel, $h$ a bandwidth, and $\hat{f}_W(\cdot)$ and $\hat{m}(\cdot)$ are defined as in (3.2) and (3.3) respectively (except that we now use the notation $\hat{m}(\cdot)$ instead of $\hat{g}(\cdot)$). Finally, the estimator $\hat{\theta}$ is defined by

$$\hat{\theta} = \hat{S}^{-1} n^{-1} \sum_{i=1}^n w^2(W_i)[Z_i - \hat{m}_Z(W_i)][Y_i - \hat{m}(W_i)],$$

where

$$\hat{S} = n^{-1} \sum_{i=1}^n w^2(W_i)[Z_i - \hat{m}_Z(W_i)][Z_i - \hat{m}_Z(W_i)]^t.$$

Define $T(x) = E[\text{Var}(Y|X) w^2(W) I(X \in J_x)]$, and consider the following assumptions :

**(PL1)** Conditions (1)–(6) given in [25], p. 72-73, are valid.

**(PL2)** $\inf_{x \in R_X} T(x) > 0$ and $\sup_{x \in R_X} T(x) < \infty$.

The conditions in [25] guarantee the weak convergence of the process $R_n^{PL}(x)$ ($x \in R_X$). They consist of certain smoothness conditions on $m(w)$, $m_Z(w)$ and $q(z|w) = P(Z \in J_z | W = w)$, conditions on the bandwidth $h$ and on the kernel $k$, moment conditions on $Y, Z$ and $\varepsilon$, and conditions on the weight function $w$.

**Theorem 3.6.** *Assume (PL1)–(PL2). Then, under $H_0^{PL}$, conditions (C0)–(C3) hold true for $V(x)$ a zero-mean Gaussian process with covariance function*

$$\text{Cov}(V(x_1), V(x_2)) = E\{\varepsilon^2 w^2(W) D(X, x_1) D(X, x_2)\},$$

*where*

$$D(X, x) = I(X \in J_x) - q(z|W) I(W \in J_w)$$
$$- E\Big\{ w(W) I(X \in J_x)[Z - m_Z(W)]^t \Big\} S^{-1} [Z - m_Z(W)],$$

*and $S = E\{w^2(W)[Z - m_Z(W)][Z - m_Z(W)]^t\}$.*

The proof of this result is very similar to the one of Theorem 3.5 and is therefore omitted.

The consistency of the bootstrap procedure described in Section 2 is established in Theorem 5.3.1 in [25]. We refer to this book for the assumptions under which this approximation is valid.

### *3.5. Other models*

The procedure developed in this paper can also be applied to a number of other situations : testing for homoscedasticity ([26]), the comparison of regression curves ([13]), single index models ([24]), autoregressive models ([18]), testing for conditional independence ([3]), among others. In each of these situations, a suitable marked empirical process has been proposed and its asymptotic behavior has been developed. In a similar way as we have done for the applications worked out in Subsections 3.1–3.4, one can use these convergence results to show the validity of the proposed empirical likelihood tests. We leave the details to the reader.

## 4. Simulations

We present some simulations carried out to assess the behavior of the new test for moderate sample sizes. We consider situations under the null, which allows to assess the accuracy of the bootstrap approximation, and under the alternative, to show the power. The new test is compared with the original test based on the marked empirical process. The simulated examples were chosen for the problem of testing a parametric regression model, and for a generalized linear model with binomial response.

### *4.1. Parametric model*

The data are taken from the following regression model

$$Y_i = X_i + d(X_i) + \sigma(X_i)\varepsilon_i \qquad i \in \{1,\ldots,n\},$$

where $n$ is the sample size, $X_i$ is a uniform random variable on the interval $[0,1]$, $\varepsilon_i$ is a standard normal random variable, the parametric model to be tested is

$$H_0 : m(x) = \theta x$$

for some parameter $\theta$ to be estimated, $d$ is a function which represents the deviation from the null hypothesis and $\sigma$ is the conditional standard deviation of $Y_i$ given $X_i$. The simulations were carried out under the null hypothesis, $d = 0$, and for different alternative functions $d$. Homoscedastic models were considered by means of a constant function $\sigma$, and several heteroscedastic models were also studied by means of different functions $\sigma$.

The following tables contain the percentage of rejections for ten thousand samples, under the nominal level 5%. The results are given for different values of the sample size, $n$, for different alternative functions $d$ and for different standard deviations of the error. The functions $d$ are coded as 0. $d(x) = 0$, 1. $d(x) = x^2$, 2. $d(x) = 0.3x\exp(x)$, 3. $d(x) = 0.3\sin(4\pi x)$ and 4. $d(x) = 0.4xI(x \leq 0.5) - 0.4(1-x)I(x > 0.5)$. The standard deviation is coded as 1. $\sigma(x) = 0.25$, 2. $\sigma(x) = 0.5x$ and 3. $\sigma(x) = 0.125(2-x)$. The columns "IRF-KS" and "IRF-CVM" contain the percentage of rejections for the test by [16], with the Kolmogorov-Smirnov and Cramér-von Mises statistics, respectively. In order to approximate the critical values for this test, wild bootstrap was used, as it is explained in [17]. The columns "EL-KS" and "EL-CVM" contain the percentage of rejections for the new test, with the Kolmogorov-Smirnov and Cramér-von Mises statistics, respectively. The sets were constructed using the value $a = 0.5$. In order to approximate the critical values, the bootstrap proposed in this paper was used. Five thousand bootstrap replicates were used for bootstrap approximations.

TABLE 1
*Percentage of rejections in homoscedastic models*

| $\sigma$ | $d$ | $n$ | IRF-KS | IRF-CVM | EL-KS | EL-CVM |
|---|---|---|---|---|---|---|
| 1 | 0 | 50 | 5.35 | 5.09 | 7.07 | 5.90 |
| 1 | 0 | 100 | 5.09 | 5.11 | 5.80 | 5.27 |
| 1 | 1 | 50 | 73.07 | 68.51 | 73.07 | 76.30 |
| 1 | 1 | 100 | 95.90 | 93.69 | 95.23 | 96.72 |
| 1 | 2 | 50 | 31.76 | 28.26 | 32.65 | 33.94 |
| 1 | 2 | 100 | 57.69 | 51.71 | 55.07 | 59.50 |
| 1 | 3 | 50 | 46.37 | 37.65 | 37.49 | 35.78 |
| 1 | 3 | 100 | 79.60 | 68.76 | 70.06 | 66.85 |
| 1 | 4 | 50 | 40.72 | 29.87 | 50.09 | 32.06 |
| 1 | 4 | 100 | 73.73 | 54.25 | 81.24 | 56.79 |

Under the null, the percentages of rejections are close to the nominal level 5% in all cases, the main deviation coming from the empirical likelihood tests with sample size $n = 50$, while for $n = 100$ the percentages are already close to 5%. Under the alternative, the power of the empirical likelihood tests is generally higher than the power of the original tests. This is the case for models 1, 2 and 4, while for model 3, a sinusoidal alternative, the power of the orignal test is higher. The improvement obtained by means of the empirical likelihood idea is larger under heteroscedastic models.

TABLE 2
*Percentage of rejections in heteroscedastic models*

| $\sigma$ | $d$ | $n$ | IRF-KS | IRF-CVM | EL-KS | EL-CVM |
|---|---|---|---|---|---|---|
| 2 | 0 | 50 | 6.11 | 5.82 | 9.01 | 7.02 |
| 2 | 0 | 100 | 5.33 | 5.17 | 6.94 | 5.43 |
| 2 | 1 | 50 | 76.00 | 83.85 | 85.67 | 83.77 |
| 2 | 1 | 100 | 96.08 | 98.52 | 98.40 | 98.14 |
| 2 | 2 | 50 | 34.57 | 39.60 | 45.69 | 42.02 |
| 2 | 2 | 100 | 58.89 | 66.76 | 67.29 | 65.82 |
| 2 | 3 | 50 | 48.13 | 47.67 | 37.01 | 37.50 |
| 2 | 3 | 100 | 79.27 | 79.76 | 68.06 | 65.10 |
| 2 | 4 | 50 | 37.48 | 40.79 | 68.25 | 41.22 |
| 2 | 4 | 100 | 72.26 | 70.67 | 95.69 | 67.92 |

TABLE 3
*Percentage of rejections in heteroscedastic models.*

| $\sigma$ | $d$ | $n$ | IRF-KS | IRF-CVM | EL-KS | EL-CVM |
|---|---|---|---|---|---|---|
| 3 | 0 | 50  | 4.92  | 4.93  | 6.84  | 5.68 |
| 3 | 0 | 100 | 4.70  | 4.51  | 5.20  | 4.86 |
| 3 | 1 | 50  | 84.88 | 79.23 | 84.54 | 87.91 |
| 3 | 1 | 100 | 99.38 | 98.19 | 99.21 | 99.64 |
| 3 | 2 | 50  | 41.28 | 36.15 | 42.22 | 46.12 |
| 3 | 2 | 100 | 71.62 | 63.10 | 70.04 | 75.65 |
| 3 | 3 | 50  | 53.50 | 42.86 | 46.42 | 43.23 |
| 3 | 3 | 100 | 88.54 | 78.65 | 88.65 | 83.21 |
| 3 | 4 | 50  | 56.09 | 38.69 | 64.88 | 43.68 |
| 3 | 4 | 100 | 87.60 | 66.19 | 92.66 | 72.83 |

## 4.2. Generalized linear model

The simulated responses $Y_i$ are taken from the binomial distribution with 15 trials and the probability of success $p(X_i)$ for $i \in \{1, \ldots, n\}$, where $n$ is the sample size and $X_i$ is a uniform random variable on the cube $[-1,1] \times [-1,1] \times [0,2]$.

The generalized linear model to be tested is

$$H_0 : p(x) = \gamma\left(\beta^t x\right)$$

where $\gamma(x) = e^x/(1+e^x)$ is the logistic transformation and $\beta$ is some parameter to be estimated.

For the simulated model under the null, the parameter $\beta$ was defined to be $\beta = (1, 2, 0.5)^t$. Two alternatives were considered: $p(x) = \Phi(\beta^t x)$ where $\Phi$ is the standard normal distribution function, and $\beta = (1, 2, 0.5)^t$ as before, and $p(x) = \gamma(x_1 + 2x_2 + 0.25(x_2+1)^2)$ with $\gamma$ the logistic transformation. We call them a probit model and a quadratic logistic model, respectively.

Table 4 contains the percentage of rejections for ten thousand samples, under the nominal level 5%. The results are given for different values of the sample size $n$, and for the three models. The five columns named "M-IRF-CVM" contain the percentages of rejections for the test by [20], which is a test based on a martingale-transformed Cramér-von Mises statistic. The parameters were estimated using the maximum likelihood method. For the transformation, smoothing is necessary in order to estimate the function $E(X|\beta^t X = u)$. To this aim, a Nadaraya-Watson estimator was used with Epanechnikov kernel on the interval $[-1, 1]$, and five possible values of the bandwidth: 1.0, 1.5, 2.0, 3.0 and 4.0. The five columns in Table 4 contain the results for these five bandwidths, respectively.

The columns "EL-KS" and "EL-CVM" contain the percentage of rejections for the test proposed in this paper, with the Kolmogorov-Smirnov and Cramérvon Mises statistics, respectively. The sets were constructed using the value $a = 0.5$. In order to approximate the critical values, the bootstrap based on multipliers, as proposed in Section 2, was used. Five thousand bootstrap replicates were used for bootstrap approximations.

TABLE 4
*Percentage of rejections in models with binomial response*

| Model | $n$ | M-IRF-CVM | | | | | EL-KS | EL-CVM |
|---|---|---|---|---|---|---|---|---|
| | | 1.0 | 1.5 | 2.0 | 3.0 | 4.0 | | |
| Null | 50 | 4.87 | 4.89 | 5.04 | 4.97 | 5.04 | 7.53 | 5.70 |
| Null | 100 | 4.19 | 4.08 | 4.15 | 4.13 | 4.03 | 7.01 | 5.31 |
| Null | 500 | 4.84 | 4.84 | 4.75 | 4.66 | 4.63 | 5.47 | 4.87 |
| Probit | 50 | 0.33 | 0.31 | 0.31 | 0.27 | 0.31 | 13.32 | 11.86 |
| Probit | 100 | 0.40 | 0.42 | 0.38 | 0.38 | 0.37 | 13.33 | 12.96 |
| Probit | 500 | 7.51 | 7.50 | 7.44 | 8.02 | 8.73 | 21.10 | 36.85 |
| Quadratic | 50 | 8.36 | 8.06 | 8.23 | 7.93 | 8.03 | 12.79 | 12.34 |
| Quadratic | 100 | 13.44 | 13.18 | 12.71 | 12.25 | 12.20 | 16.44 | 19.84 |
| Quadratic | 500 | 68.05 | 66.60 | 64.44 | 62.31 | 62.34 | 49.53 | 75.69 |

These results show a good approximation to the level under the null. The biggest deviations come from the Kolmogorov-Smirnov statistic of the empirical likelihood process with small sample size, although the result for $n = 500$ shows the consistency of the bootstrap. The power is much larger for the empirical likelihood test, especially with the Cramér-von Mises statistic. Note that the martingale-transformed statistic is based on a Cramér-von Mises statistic. The difficulty to detect the probit alternative should also be emphasized, with percentages of rejections even smaller than the nominal level for the martingale-transformed statistic. By other way, the bandwidth has not much influence on the results for the martingale-transformed statistic.

## Appendix: Proofs

In this Appendix we prove the main result (Theorem 2.1) that states the generic conditions (C0)–(C3) under which the convergence of the test statistics is guaranteed, and we also verify these conditions for the models considered in Section 3.

**Proof of Theorem 2.1.**
(a) The result is an easy extension of Theorem 2.1 in [11], with $a_n = 1$ and $h = (\theta, g)$. The extension lies in the fact that their result is limited to the convergence of $\ell(u, \hat\theta, \hat g)$ for a fixed value of $u$, whereas we consider a process in $u$. It is however a straightforward exercise to show that under conditions (C0)–(C3), the process $\ell(u, \hat\theta, \hat g)$ ($u \in U$) converges to $W(u)$ in $\ell^\infty(U)$. The convergence of the test statistics $S_n$ and $T_n$ now follows immediately.

(b) We will show that $P(S_n > c) \to 1$ as $n$ tends to infinity. The proof for $T_n$ is similar. Write $P(S_n > c) = P(c < S_n < \infty) + P(S_n = \infty)$. We will show that $P(c < S_n < \infty) - P(S_n < \infty) \to 0$, from which it then follows that $P(S_n > c) \to 1$. If $S_n < \infty$, then for each $u \in U$, the supremum in (2.4) is a supremum over a non-empty set, and hence it is easy to show that

$$\ell(u, \hat\theta, \hat g) = 2\sum_{i=1}^n \log\Big(1 + \lambda_n(u) A_i(u)\Big),$$

where the Lagrange multiplier $\lambda_n(u)$ satisfies

$$\sum_{i=1}^{n} \frac{A_i(u)}{1 + \lambda_n(u) A_i(u)} = 0, \tag{A.1}$$

and $A_i(u) = L(X_i, u, \hat{\theta}, \hat{g})[Y_i - \gamma(X_i, \hat{\theta}, \hat{g})]$ $(i = 1, \ldots, n)$. First, we will show that under $H_A$, and provided $S_n < \infty$,

$$n^{-1/2} \Big( \inf_{u \in U_0} |\lambda_n(u)| \Big)^{-1} = o_P(1),$$

i.e. for all $\eta > 0$,

$$P\Big(n^{1/2} \inf_{u \in U_0} |\lambda_n(u)| > \eta \,\Big|\, S_n < \infty \Big) \to 1, \tag{A.2}$$

where $U_0$ is the subset of $U$ defined in condition (C1'). Note that it follows from (A.1) that

$$0 = \sum_{i=1}^{n} A_i(u) - \lambda_n(u) \sum_{i=1}^{n} \frac{A_i^2(u)}{1 + \lambda_n(u) A_i(u)},$$

and hence

$$n^{1/2} \inf_{u \in U_0} |\lambda_n(u)| = n^{1/2} \inf_{u \in U_0} \Big| \frac{n^{-1} \sum_{i=1}^{n} A_i(u)}{n^{-1} \sum_{i=1}^{n} [A_i^2(u)][1 + \lambda_n(u) A_i(u)]^{-1}} \Big|$$
$$\geq \frac{\inf_{u \in U_0} |n^{-1} \sum_{i=1}^{n} A_i(u)|}{n^{-1/2} \sup_{u,i} |A_i(u)|^2},$$

since $n^{-1} \sum_{i=1}^{n} [1 + \lambda_n(u) A_i(u)]^{-1} = 1$, and this goes to infinity, using (C1') and (C3'). Hence, (A.2) holds true. Next write

$$\ell(u, \hat{\theta}, \hat{g}) = 2\lambda_n(u) \sum_{i=1}^{n} \frac{A_i(u)}{1 + \lambda_n(u) A_i(u)} + \lambda_n^2(u) \sum_{i=1}^{n} \frac{A_i^2(u)}{(1 + \xi_i(u))^2},$$

for some $|\xi_i(u)| \leq |\lambda_n(u)||A_i(u)|$. Note that the first term on the right hand side equals 0. Hence,

$$\sup_{u \in U_0} \ell(u, \hat{\theta}, \hat{g}) \geq \frac{\sup_{u \in U_0} n^{-1} \sum_{i=1}^{n} A_i^2(u)}{2 \sup_{u \in U_0} (|n^{1/2} \lambda_n(u)|^{-2}) + 2 \sup_{u,i} |n^{-1/2} A_i(u)|^2}.$$

It now follows that

$$P(c < S_n < \infty) \geq P(\sup_{u \in U_0} \ell(u, \hat{\theta}, \hat{g}) > c | S_n < \infty) P(S_n < \infty)$$
$$\geq P\Big( \sup_u |T(u)| - R_{n1} > 2c \sup_{u \in U_0} (|n^{1/2} \lambda_n(u)|^{-2}) + 2c R_{n2},$$
$$R_{n1} \leq \varepsilon_1, R_{n2} \leq \varepsilon_2 \,\Big|\, S_n < \infty \Big) P(S_n < \infty), \tag{A.3}$$

for some $\varepsilon_1, \varepsilon_2 > 0$ to be chosen later, where $R_{n1} = \sup_u |n^{-1} \sum_{i=1}^n A_i^2(u) - T(u)|$ and $R_{n2} = \sup_{u,i} |n^{-1/2} A_i(u)|^2$. Next, (A.3) is bounded below by

$$P\Big( \inf_{u \in U_0} |n^{1/2} \lambda_n(u)| > (2c)^{1/2} \Big[ \sup_{u \in U} |T(u)| - \varepsilon_1 - 2c\varepsilon_2 \Big]^{-1/2},$$

$$R_{n1} \leq \varepsilon_1, R_{n2} \leq \varepsilon_2 \Big| S_n < \infty \Big) P(S_n < \infty)$$

$$\geq P\Big( \inf_{u \in U_0} |n^{1/2} \lambda_n(u)| > (2c)^{1/2} \Big[ \sup_{u \in U} |T(u)| - \varepsilon_1 - 2c\varepsilon_2 \Big]^{-1/2} \Big| S_n < \infty \Big) P(S_n < \infty)$$

$$+ P(R_{n1} \leq \varepsilon_1) + P(R_{n2} \leq \varepsilon_2) - 2$$

$$= P(S_n < \infty) + o(1),$$

for $\varepsilon_1, \varepsilon_2 > 0$ small enough, using (C2), (C3') and (A.2). This finishes the proof.

**Proof of Theorem 3.1.**

For condition (C0), note that

$$P(EL(x, \hat{\theta}) = 0 \text{ for some } x \in R_X)$$

$$= P(\text{for some } x \in R_X, Y_i - \gamma(X_i, \hat{\theta}) > 0 \text{ for all } X_i \in J_x$$

$$\text{or } Y_i - \gamma(X_i, \hat{\theta}) < 0 \text{ for all } X_i \in J_x)$$

$$\leq P(\varepsilon_i > -\delta \text{ for all } X_i \leq a \text{ or for all } X_i \geq a)$$

$$+ P(\varepsilon_i < \delta \text{ for all } X_i \leq a \text{ or for all } X_i \geq a)$$

$$+ P(\sup_{x \in R_X} |\gamma(x, \hat{\theta}) - \gamma(x, \theta_0)| > \delta)$$

$$\leq p_1^{l_n} + p_2^{l_n} + p_1^{k_n} + p_2^{k_n} + o(1)$$

$$= o(1),$$

where $\delta > 0$, $p_1 = P(\varepsilon > -\delta)$, $p_2 = P(\varepsilon < \delta)$, $l_n = \sum_{i=1}^n I(X_i \leq a)$ and $k_n = n - l_n$. The latter equality follows from the fact that $p_1 < 1$ and $p_2 < 1$ for $\delta > 0$ small enough, and the fact that $k_n, l_n \to \infty$ since $a$ belongs to the interior of $R_X$. Condition (C1) follows from Corollary 1.3 in [16] (it is easily seen how $I(X_i \leq x)$ can be replaced by $I(X_i \in J_x)$ in the proof of that corollary $(i = 1, \ldots, n)$). For (C2), consider

$$n^{-1} \sum_{i=1}^n I(X_i \in J_x)[Y_i - \gamma(X_i, \hat{\theta})]^2$$

$$= n^{-1} \sum_{i=1}^n I(X_i \in J_x)[Y_i - \gamma(X_i, \theta_0)]^2$$

$$+ n^{-1} \sum_{i=1}^n I(X_i \in J_x) \Big( \frac{\partial \gamma(X_i, \hat{\theta}_1)}{\partial \theta} \Big)^2 (\hat{\theta} - \theta_0)^2$$

$$- 2n^{-1} \sum_{i=1}^n I(X_i \in J_x)[Y_i - \gamma(X_i, \theta_0)] \frac{\partial \gamma(X_i, \hat{\theta}_1)}{\partial \theta} (\hat{\theta} - \theta_0)$$

$$= n^{-1} \sum_{i=1}^n I(X_i \in J_x)[Y_i - \gamma(X_i, \theta_0)]^2 + o_P(1), \quad (A.4)$$

uniformly in $x \in R_X$, provided conditions (P1)–(P3) hold, and where $\hat{\theta}_1$ is between $\theta_0$ and $\hat{\theta}$. To show that (A.4) converges to $T(x)$ uniformly in $x$, it suffices to show that the class of functions

$$\mathcal{F} = \{(t, y) \to I(t \in J_x)[y - \gamma(t, \theta_0)]^2 : x \in R_X\}$$

is Glivenko-Cantelli. This is obvious, since the class of functions $t \to I(t \in J_x)$ is contained in the class of monotone and bounded functions into $[0, 1]$, which is Glivenko-Cantelli (see [22], p. 149), and since the function $(t, y) \to [y - \gamma(t, \theta_0)]^2$ is independent of $x$. The fact that $\inf_{x \in R_X} T(x) > 0$ and that $\sup_{x \in R_X} T(x) < \infty$ follows from condition (P1). We finally prove (C3). Note that the left hand side of the statement in condition (C3) is bounded by

$$\max_{1 \leq i \leq n} |Y_i - \gamma(X_i, \theta_0)| + \max_{1 \leq i \leq n} |\gamma(X_i, \hat{\theta}) - \gamma(X_i, \theta_0)|.$$

The first expression above is $o_P(n^{1/2})$ (see Lemma 11.2 in [14]), whereas the second term is $o_P(1)$ by the uniform boundedness of $\frac{\partial \gamma(x, \theta)}{\partial \theta}$ over $x$ and $\theta$ and since $\hat{\theta} - \theta_0 = o_P(1)$.

**Proof of Theorem 3.2.**
We only need to show the validity of condition (C4). Write

$$R_n^{P*}(x) - \tilde{R}_n^{P*}(x) = -n^{-1/2} \sum_{i=1}^n I(X_i \in J_x)[\gamma(X_i, \hat{\theta}) - \gamma(X_i, \theta_0)]V_i$$
$$- G^t(x, \theta_0) n^{-1/2} \sum_{i=1}^n [\hat{h}(X_i, Y_i, \hat{\theta}) - h(X_i, Y_i, \theta_0)]V_i$$
$$- [\hat{G}^t(x, \hat{\theta}) - G^t(x, \theta_0)] n^{-1/2} \sum_{i=1}^n h(X_i, Y_i, \theta_0)V_i + o_P(1).$$

The second term above is $o_P(1)$ by condition $(P4)$. The convergence to zero of the first term follows from the fact that by assumption $(P5)$,

$$n^{-1/2} \sum_{i=1}^n I(X_i \in J_x)[\gamma(X_i, \hat{\theta}) - \gamma(X_i, \theta_0)]V_i$$
$$= (\hat{\theta} - \theta_0)^t n^{-1/2} \sum_{i=1}^n I(X_i \in J_x) \frac{\partial}{\partial \theta} \gamma(X_i, \theta_0)V_i + o_P(1) = o_P(1),$$

since $E[I(X_i \in J_x) \frac{\partial}{\partial \theta} \gamma(X_i, \theta_0) V_i] = 0$. Finally, for the third term above note that $n^{-1/2} \sum_{i=1}^n h(X_i, Y_i, \theta_0) V_i = O_P(1)$ and that $\hat{G}(x, \hat{\theta}) - G(x, \theta_0) = \hat{G}(x, \theta_0) - G(x, \theta_0) + o_P(1)$. In a similar way as in the proof of Theorem 3.1 one can show that the class

$$\left\{t \to I(t \in J_x) \frac{\partial}{\partial \theta} \gamma(t, \theta_0) : x \in R_X\right\}$$

is Glivenko-Cantelli. Hence, the third term is also $o_P(1)$.

**Proof of Theorem 3.3.**
Condition (C0) can be shown in a similar way as in the proof of Theorem 3.1. Condition (C1) is an immediate consequence of Theorem 1 in [20]. For (C2), note that

$$n^{-1} \sum_{i=1}^{n} I(\hat{\beta}^t X_i \leq x)[Y_i - \gamma(\hat{\beta}^t X_i, \hat{\alpha})]^2$$

$$= n^{-1} \sum_{i=1}^{n} I(\beta_0^t X_i \leq x)[Y_i - \gamma(\beta_0^t X_i, \alpha_0)]^2$$

$$+ n^{-1} \sum_{i=1}^{n} I(\beta_0^t X_i \leq x) \left[ \frac{\partial \gamma(\beta_0^t X_i, \hat{\alpha}_1)}{\partial \alpha} \right]^2 (\hat{\alpha} - \alpha_0)^2$$

$$- 2n^{-1} \sum_{i=1}^{n} I(\beta_0^t X_i \leq x)[Y_i - \gamma(\beta_0^t X_i, \alpha_0)] \frac{\partial \gamma(\beta_0^t X_i, \hat{\alpha}_1)}{\partial \alpha} (\hat{\alpha} - \alpha_0)$$

$$+ n^{-1} \sum_{i=1}^{n} [I(\hat{\beta}^t X_i \leq x) - I(\beta_0^t X_i \leq x)][Y_i - \gamma(\hat{\beta}^t X_i, \hat{\alpha})]^2$$

$$+ n^{-1} \sum_{i=1}^{n} I(\beta_0^t X_i \leq x) \{[Y_i - \gamma(\hat{\beta}^t X_i, \hat{\alpha})]^2 - [Y_i - \gamma(\beta_0^t X_i, \hat{\alpha})]^2\},$$

where $\hat{\alpha}_1$ is between $\alpha_0$ and $\hat{\alpha}$. The first three terms above can be dealt with in a similar way as in the proof of Theorem 3.1. For the fourth term above, note that, since $\hat{\beta} - \beta_0 = O_P(n^{-1/2})$ (by condition (GLM2)), the number of $i$-values for which $I(\hat{\beta}^t X_i \leq x) - I(\beta_0^t X_i \leq x) \neq 0$ is of the order $O((n \log \log n)^{1/2})$ a.s. (by the law of iterated logarithm, see e.g. [15], p. 35). Hence, the fourth term is $O((n^{-1} \log \log n)^{1/2})$ a.s. The fifth term is $O_P(n^{-1/2})$, by conditions (GLM1), (GLM2) and (GLM3)(ii). Finally, to check (C3), a similar derivation as in the proof of Theorem 3.1 can be followed.

**Proof of Theorem 3.4.**
The proof can be established in a somewhat analogous way as for the parametric model. We therefore restrict attention to the main points of difference. In a similar way as in the proof of Theorem 1 in [20], it can be easily shown that when one replaces the indicator $I(\hat{\beta}^t X_i \in J_x)$ by $I(\beta_0^t X_i \in J_x)$ in the formula of $R_n^{GLM*}(u)$ (i.e. in the formulas of $\hat{Q}(X_i, Y_i, u, \hat{\theta})$ and $\hat{G}(u, \hat{\theta})$), the asymptotic distribution does not change. This so-obtained process can now be analyzed in a very similar way as the process $R_n^{P*}(x)$ in the proof of Theorem 3.2, which then shows the equivalence of $R_n^{GLM*}(u)$ and $\tilde{R}_n^{GLM*}(u)$, uniformly in $u$.

**Proof of Theorem 3.5.**
For condition (C0), a derivation similar to the one in the proof of Theorem 3.1 can be given, using the fact that $\sup_w |\hat{g}(w) - g(w)| = o_P(1)$. Condition (C1) follows from Theorem 1 in [3]. For condition (C2), write

$$n^{-1} \sum_{i=1}^{n} I(X_i \in J_x) \hat{f}_W^2(W_i)(Y_i - \hat{g}(W_i))^2$$
$$= n^{-1} \sum_{i=1}^{n} I(X_i \in J_x) f_W^2(W_i)(Y_i - g(W_i))^2$$
$$+ n^{-1} \sum_{i=1}^{n} I(X_i \in J_x)[\hat{f}_W^2(W_i) - f_W^2(W_i)](Y_i - g(W_i))^2$$
$$+ n^{-1} \sum_{i=1}^{n} I(X_i \in J_x) \hat{f}_W^2(W_i)[(Y_i - \hat{g}(W_i))^2 - (Y_i - g(W_i))^2].$$

The first term above can be treated as in the proof of Theorem 3.1. The second one is $o_P(1)$ by using assumption $(SV4)$ and the fact that $\sup_w |\hat{f}_W(w) - f_W(w)| = o_P(1)$, while the negligeability of the third term follows again from condition $(SV4)$ and from the uniform consistency of $\hat{g}(w)$ for all $w \in R_W$. Finally, for (C3), note that

$$\sup_x \max_i |I(X_i \in J_x) \hat{f}_W(W_i)(Y_i - \hat{g}(W_i))|$$
$$\leq \max_i |f_W(W_i)(Y_i - g(W_i))| + \max_i |[\hat{f}_W(W_i) - f_W(W_i)](Y_i - \hat{g}(W_i))|$$
$$+ \max_i |f_W(W_i)[\hat{g}(W_i) - g(W_i)]|$$
$$= o_P(1),$$

by using again the same development as in the proof of Theorem 3.1 for the first term, and the uniform consistency of $\hat{g}$ and $\hat{f}_W$ for the second and third term.